\documentclass[12pt,oneside,english]{amsart}
\usepackage[latin1]{inputenc}
\usepackage{amssymb}

\makeatletter
\newtheorem{theorem}{Theorem}
\newtheorem{lemma}{Lemma}
\newtheorem{example}{Example}
\newtheorem{definition}{Definition}
\newtheorem{remark}{Remark}

\newtheorem{corollary}{Corollary}

\usepackage{babel}

\makeatother
\begin{document}
\baselineskip=17pt

\title[ Sequences which are similar to the positive integers]{Several results on sequences which are similar to the positive integers}

\author{Vladimir Shevelev}
\address{Departments of Mathematics \\Ben-Gurion University of the
 Negev\\Beer-Sheva 84105, Israel. e-mail:shevelev@bgu.ac.il}


\begin{abstract}
Sequence of positive integers $\{x_n\}_{n\geq1}$ is called similar to $\mathbb {N}$ with respect to a given property $A$ if for every $n\geq1$ the numbers $x_n$ and $n$ are in the same class of equivalence with respect to $A\enskip(x_n\sim n (prop \enskip A)).$ If $x_1=a(>1)\sim1 (prop\enskip A)$ and $x_n>x_{n-1}$ with the condition that $x_n$ is the nearest to $x_{n-1}$ number such that $x_n\sim n (prop \enskip A),$ then the sequence $\{x_n\}$ is called minimal recursive with the first term $a\enskip(\{x_n^{(a)}\}).$ We study two cases: $A=A_1$ is the value of exponent of the highest power of 2 dividing an integer and $A=A_2$ is the parity of the number of ones in the binary expansion of an integer. In the first case we prove that, for sufficiently large $n, \enskip x_n^{(a)}=x_n^{(3)};$ in the second case we prove that, for $a>4$ and sufficiently large $n,\enskip x_n^{(a)}=x_n^{(4)}.$

\end{abstract}

\maketitle

\section{Introduction, main notions and results}
   Two positive integers $x, y$ which are in the same class of equivalence with respect to a given property $A$ are called similar respect to $A,$ denoting this by $x\sim y \enskip (prop\enskip  A).$ Two sequences of positive integers $\{x_n\}_{n\geq1}, \{y_n\}_{n\geq1}$ are called similar with respect to$A$ ($\{x_n\}_{n\geq1}\sim \{y_n\}_{n\geq1}\enskip (prop\enskip  A)$) if $x_n\sim y_n (prop A)$ for every $n=1,2,...$ Let, e.g.,$A=A_1$ be the value of exponent of the highest power of 2 dividing an integer and $A=A_2$ be the parity of the number of ones in the binary expansion of an integer. Well known Moser-de Bruijn sequence, that is ordered sums of distinct powers of 4 (see A000695 in [2]), gives an example of a similar to $\mathbb {N}$ sequence with respect to property $A_2:$
\begin{equation}\label{1}
1, 4, 5, 16, 17, 20, 21, 64, 65, 68, 69, 80, 81, 84, 85, 256, 257, 260, 261, ...
\end{equation}
A non-trivial example of a fast growing similar to $\mathbb {N}$ sequence with respect to property $A_1$ is given by the following theorem.
\begin{theorem}\label{1}For $n\geq1,$
\begin{equation}\label{2}
((2n-1)!!+(-1)^{(n-1)(n-2)/2})/2 \sim\mathbb {N} \enskip(prop \enskip A_1).
\end{equation}
\end{theorem}
It is the sequence
\begin{equation}\label{3}
1, 2, 7, 52, 473, 5198, 67567, 1013512, 17229713, 327364538, ...
\end{equation}
\begin{remark}
Theorem 1 is a part of our full research of the binary carry sequence ($A007814$ in $[2])$ of $(2n-1)!!\pm1$ (see explicit formulas to our sequences $A158570, A158572$ in $[2]$).
 \end{remark}
 Let $x_1=a\sim 1 \enskip (prop\enskip  A)$ and $x_n>x_{n-1}$ with the condition that $x_n$ is the nearest to $x_{n-1}$ number such that $x_n\sim n \enskip (prop\enskip  A).$ Then the sequence $\{x_n\}$ is called minimal recursive similar to $\mathbb {N}$ with the first term $a$ and we denote it by $(\{x_n^{(A;\enskip a)}\}).$ Finally, two sequences $\{x_n\}$ and $\{y_n\}$ are called essentially coincide if, for all sufficiently large $n,$ we have $x_n=y_n.$ Consider some examples.
\begin{example}\label{1}
Let $A=A_1$ be the value of exponent of the highest power of 2 dividing an integer.
\end{example}
Then $3\sim 1 \enskip (prop\enskip  A_1)$ and the first terms of $\{x_n^{(A_1;\enskip3)}\}$ are:
    \begin{equation}\label{4}
 3,6,7,12,13,14,15,24,25,26,27,28,29,30,31,48,49,50,51,52,...
\end{equation}
The first main our result is following.
\begin{theorem}\label{2}
Every minimal recursive sequence $\{x_n^{(A_1;\enskip a)}\}$ essentially coincides with sequence $\{x_n^{(A_1;\enskip3)}\}$ $(4).$
\end{theorem}
\begin{example}\label{2}
Let $A=A_2$ be the parity of the number of ones in the binary expansion of an integer.
\end{example}
Then $2\sim 1 \enskip (prop \enskip A_2)$ and the first terms of $\{x_n^{(A_2;\enskip2)}\}$ are:
\begin{equation}\label{5}
2,4,5,7,9,10,11,13,15,17,19,20,21,22,23,25,27,29,31,33,...
\end{equation}
Analogously,
$4\sim 1 \enskip (prop \enskip A_2)$ and the first terms of $\{x_n^{(A_2;\enskip4)}\}$ are:
\begin{equation}\label{6}
4,7,9,11,12,15,16,19,20,23,25,27,28,31,33,35,36,39,41,43,...
\end{equation}
The second main our result is following.
\begin{theorem}\label{3}
Every minimal recursive sequence $\{x_n^{(A_2;\enskip a)}\}$ either is sequence $\{x_n^{(A_2;\enskip 2)}\}$ $(5)$ or essentially coincides with sequence $\{x_n^{(A_2;\enskip 4)}\}$ $(6).$
\end{theorem}
Evidently, the number of examples could be continued infinitely. We give only two more.
\begin{example}\label{3}
Let $A=A_3$ be the property of a number to be or not to be prime.
\end{example}
Then  $4\sim 1 \enskip (prop \enskip A_3).$ The first terms of the corresponding minimal recursive sequence $\{x_n^{(A_3;\enskip4)}\}$ are:
\begin{equation}\label{7}
4,5,7,8,11,12,13,14,15,16,17,18,19,20,21,22,23,24,29,30,32,...
\end{equation}\newpage
It is interesting that here the sequences $ \{x_n^{(A_3;\enskip6)}\}, \{x_n^{(A_3;\enskip10)}\},\{x_n^{A_3;\enskip(12)}\}$ essentially coincide with $\{x_n^{A_3;\enskip(4)}\},$ but the question "whether the sequence $\{x_n^{A_3;\enskip(16)}\}$ essentially coincides with $\{x_n^{A_3;\enskip(4)}\}$?" remains open.
\begin{example}\label{4}
Let $A=A_4$ be the number of prime divisors of an integer.
\end{example}
Then $3 \sim 2 \enskip (prop \enskip A_4).$ Here there is a sense to consider minimal recursive similar to $\mathbb {N}\setminus\{1\}\enskip(prop\enskip  A_4)$ with the first term $a,$ and we denote it by $(\{x_n^{(A_4;\enskip 2;\enskip a)}\}).$ The first terms of the corresponding minimal recursive sequence $\{x_n^{(A_4;\enskip 2;\enskip3)}\}$are:
\begin{equation}\label{8}
3,5,7,8,10,11,13,16,18,19,20,23,24,26,27,29,33,37,38,...
\end{equation}
We verified that the sequences $$\{x_n^{(A_4;\enskip 2;\enskip 5)}\}, \{x_n^{(A_4;\enskip 2;\enskip 7)}\}, \{x_n^{(A_4;\enskip 2;\enskip 8)}\},\{x_n^{(A_4;\enskip 2;\enskip 9)}\}$$ essentially coincide with $\{x_n^{(A_4;\enskip 2;\enskip  3)}\},$ but we do not know whether $\{x_n^{(A_4;\enskip 2;\enskip a)}\}$ essentially coincides with $\{x_n^{(A_4;\enskip 2;\enskip 3)}\}$ for every $a\sim3\enskip (prop \enskip A_4).$\newline
In connection with Theorem 2,3 and examples 3,4, it is natural to pose the following general problem.\newline
\bfseries Problem. \enskip\mdseries To find a characterization of the class of properties $A$ such that, for a given $n_0,$ there exists $t=t(A,n_0)\sim n_0(\geq1)\enskip(prop \enskip A)$ such that for every $a\geq t, \enskip a\sim n_0\enskip(prop \enskip A)$ the minimal recursive sequence $\{x_n^{(A;\enskip n_0; \enskip a)}\},$ similar to $\{n_0, n_0+1, ...\}$ with respect to $A,$ essentially coincides with $\{x_n^{(A;\enskip n_0;\enskip  t)}\}.$\newline
In the following sections we give proofs of Theorems 1-3. In the limits of Sections 3,4 we write $\{x_n^{(a)}\}$ for
$\{x_n^{(A_1;\enskip a)}\},$ while in the limits of Sections 5-7 we write $\{x_n^{(a)}\}$ for
$\{x_n^{(A_2;\enskip a)}\}.$
\section{Proof of Theorem 1}
Below we denote the exponent of the highest power of 2 dividing $n$ by $(n)_2.$
Using induction, distinguish the following cases: $n\equiv i \pmod 4,\enskip i=0,1,2,3.$ Note that, in cases of $i=1,2,3$ the proofs are quite analogous to the following subcase of the case $i=0:\enskip n\equiv 4 \pmod 8.$
Therefore, we prove only the case $n\equiv 0 \pmod 4$ and start with the mentioned subcase.\newline
1)Let $n=8k-4,\enskip k\geq1.$ Then $(n)_2=2$ and, according to (2), we should prove that
\begin{equation}\label{9}
((16k-9)!!-1)_2=3,\enskip k\geq1.
\end{equation}
 Denoting
$$ a_k=(16k-9)!!-1,$$
we have $(a_1)_2=(7!!-1)_2=(104)_2=3.$
Suppose that (9) is valid for some \newpage$k\geq1.$ This means that $a_k=8l,$ where $l$ is an odd number. Putting $16k=k_1,$ we have
$$a_{k+1}=(k_1+7)!!-1=$$
$$(8l+1)(k_1-7)(k_1-5)(k_1-3)(k_1-1)(k_1+1)(k_1+3)(k_1+5)(k_1+7)-1,$$
where $8l+1=(k_1-9)!!.$\newline
Consequently, $a_{k+1}$ has the form
$$a_{k+1}=16m+(8l+1)(3\cdot5\cdot7)^2-1=16m+8l\cdot105^2+105^2-1=16m+8l\cdot105^2+13\cdot53\cdot16$$
and since $l$ is odd then $(a_{k+1})_2=3.$
2)Let now $n=2^{t-1}u,\enskip t\geq4,$ where $u$ is odd. Combining this case with the previous one, we prove that
for $t\geq3,$
\begin{equation}\label{10}
((2n-1)!!-1)_2=((2^tu-1)!!-1)_2=(2n)_2=t.
\end{equation}
As the base of induction we take the case 1) which corresponds to $t=3$ and $u=2k-1, \enskip k\geq1,$ i.e. to the proved formula (9). Let (10) is true for some $t\geq3$ and every odd $u.$ Then, denoting $c_t=(2^tu-1)!!-1,$ we have
$c_t=2^tv,$ where $v$ is odd. Now we find
$$c_{t+1}=(2^{t+1}u-1)!!-1=$$
\begin{equation}\label{11}
(2^tv+1)(2^tu+1)(2^tu+3)\cdot...\cdot(2^tu+(2^tu-1))-1,
\end{equation}
where
$$2^tv+1=(2^tu-1)!!=c_t+1.$$
Choosing the first summands in at least two brackets of (11), we obtain the number of the form $2^{2t}r$ with an integer $r.$ Choosing the second summands in every bracket, beginning with the second one, we find (together with
the subtracting 1)
$$(c_t+1)(2^tu-1)!!-1=(c_t+1)^2-1=c_t^2+2c_t=2^{2t}v^2+2^{t+1}v. $$
Finally, choosing $2^tu$ in consecutive order exactly in only brackets, beginning with the second one, while in others
choosing the second summands, we obtain the following sum:
$$2^tu(2^tu-1)!!(1+1/3+1/5+...+1/(2^tu-1)=$$
$$2^tu(2^tu-1)!!\sum_{s=1}^{2^{t-1}u-1}(1/(2s-1)+1/(2s+1))=$$
$$2^{t+2}u(2^tu-1)!!\sum_{s=1}^{2^{t-1}u-1}\frac {s} {(2s-1)(2s+1)}=2^{t+2}h,$$
where $h$ is integer.
As a result, we have
$$c_{t+1}=2^{2t}(r+v^2)+2^{t+2}h+2^{t+1}v$$\newpage
with odd $v$ and thus $(c_{t+1})_2=t+1.\blacksquare$
\begin{corollary}
For every positive odd $x$ we have
\begin{equation}\label{12}
((2n-1)!!^x+(-1)^{(n-1)(n-2)/2})/2 \sim\mathbb {N} \enskip(prop \enskip A_1).
\end{equation}
\end{corollary}
\bfseries Proof. \enskip\mdseries Indeed, denoting $(2n-1)!!=a(n)$ and $(-1)^{(n-1)(n-2)/2}=b(n)=b(n)^x,$ we have
$$a(n)^x+b(n)^x=(a(n)+b(n))(a(n)^{x-1}-a(n)^{x-2}b(n)+a(n)^{x-3}b(n)^2-...+b(n)^{x-1}).$$
Since the second brackets contain odd number of odd summands, then
$$(a(n)^x+b(n)^x)_2=(a(n)^x+b(n))_2=(a(n)+b(n))_2.$$
According to Theorem 1, $(a(n)+b(n))_2=(2n)_2,$ therefore, also $(a(n)^x+b(n))_2=(2n)_2,$ and the corollary follows.$\blacksquare$
\section{Proof of Theorem 2}
\begin{lemma}\label{1}
\begin{equation}\label{13}
x_{2^t}^{(3)}=3\cdot2^t.
\end{equation}
\end{lemma}
\bfseries Proof. \enskip\mdseries
 Noting that $x_1^{(3)}=3, x_2^{(3)}=6, $ suppose that for some $t$ we have $x_{2^t}^{(3)}=3\cdot2^t.$ Then, by the definition of $\{x_n^{(3)}\}$, we have
 $$x_{2^t+j}^{(3)}=2^{t+1}+2^t+j, \enskip 0\leq j\leq 2^t-1,$$
 such that
 $$x_{2^t+2^{t}-1}^{(3)}=2^{t+1}+2^t+2^t-1.$$
 Now adding 1 to argument of $x,$  we obtain $2^{t+1},$ while, adding 1 to the right hand side, we obtain $2^{t+2}$ and, according to the algorithm, we should add $2^{t+1}.$  Thus we conclude that $x_{2^{t+1}}^{(3)}=3\cdot2^{t+1}.$ $\blacksquare$
  \begin{lemma}\label{2}
If for some $a,$ we have $ x_{2^r}^{(a)}=2^k+2^r,$ where  $1\leq r< k ,$  then there exists $T$
such that $x_{2^T}^{(a)}=3\cdot2^{T}.$
  \end{lemma}
 \bfseries Proof. \enskip\mdseries
 Let, first, $l=1.$ If $k=r+1,$ then we can take $T=r.$ Therefore, suppose that $k>r+1.$
  By the condition, $ x_{2^r}^{(a)}=2^k+2^r.$
 Evidently, we have
$$x_{2^r+j}^{(a)}=2^k+j, \enskip 0\leq j\leq 2^k-1-2^r,$$
 such that
 $$ x_{2^k-1}^{(a)}=2^{k+1}-1.$$
 Therefore, according to the algorithm of the minimal recursive sequence, we find
 $$ x_{2^k}^{(a)}=2^{k+1}+2^k$$
 and we can take $T=k.$
 $\blacksquare$\newpage
\begin{lemma}\label{3}
If for some $a,$ we have
$$ x_{2^{r_1}}^{(a)}=2^{r_t}+2^{r_{t-1}}+...+2^{r_1},\enskip r_t>r_{t-1}>...>r_1,$$
 then there exists $T=T(r_1,r_2,...,r_t)$ such that
 \begin{equation}\label{14}
 x_{2^T}^{(a)}=3\cdot2^{T}.
\end{equation}
  \end{lemma}
  \bfseries Proof. \enskip\mdseries We use induction over $t\geq2.$  The base of induction is given by Lemma 2. Suppose, that the statement is true for every $t\leq k.$ Let $t=k+1,$ such that
 $$ x_{2^{r_1}}^{(a)}=2^{r_{k+1}}+...+2^{r_3}+2^{r_2}+2^{r_1},\enskip r_{k+1}>r_{k}>...>r_1,$$
 By the minimal recursive algorithm, we have
 $$ x_{2^{r_1}+2^{r_1}-1}^{(a)}=2^{r_{k+1}}+...+2^{r_3}+2^{r_2}+2^{r_1}+2^{r_1}-1$$
 and thus
 $$ x_{2^{r_1+1}}^{(a)}=2^{r_{k+1}}+...+2^{r_3}+2^{r_2}+2^{r_1+1}.$$
 If here $r_1+1=r_2,$ then in the right hand side we have $k$ binary ones and lemma follows from the supposition.
 Suppose that $r_1+1<r_2.$ Then we find consecutively
 $$ x_{2^{r_2}-1}^{(a)}=2^{r_{k+1}}+...+2^{r_3}+2^{r_2}+2^{r_2}-1;$$
 $$ x_{2^{r_2}}^{(a)}=2^{r_{k+1}}+...+2^{r_3}+2^{r_2+1}+2^{r_2};$$
 $$ x_{2^{r_2+1}}^{(a)}=2^{r_{k+1}}+...+2^{r_3}+2^{r_2+2}+2^{r_2+1}.$$
 Note that, the both of cases $r_2+1=r_3$ and $r_2+2=r_3$ lead on the previous step or on the last one to $t=k$ and the lemma follows. Suppose that  $r_2+2<r_3.$ Finally, we obtain
  $$ x_{2^{r_2+1}+2^{r_2+1}-1}^{(a)}=2^{r_{k+1}}+...+2^{r_3}+2^{r_2+2}+2^{r_2+1}+2^{r_2+1}-1$$
  and
  $$ x_{2^{r_2+2}}^{(a)}=2^{r_{k+1}}+...+2^{r_3}+2^{r_2+3}.$$
  Now $t\leq k,$ and the lemma follows. $ \blacksquare$\newline
In particular, if $r_1=1$ and the binary expansion of $a$ has the form $a=2^{r_t}+2^{r_{t-1}}+...+2^{r_2}+2,$ then (15) is valid for some $T=T(a).$ Therefore, by Lemma 1, $$x_{2^T}^{(a)}=x_{2^T}^{(3)}$$
and, according to the minimal recursive algorithm, for $n\geq2^T,$ we have $x_n^{(a)}=x_n^{(3)}.\enskip \blacksquare$
\begin{remark}
Actually, we proved some more: if to begin the minimal recursive algorithm with an arbitrary number $n_0,$ putting
$y(n_0)=y_0>n_0,$ such that $y_0\sim n_0 \enskip (prop\enskip A_1),$ then we obtain a sequence $\{y(n)\}$ which is essentially coincides with $\{x_n^{(3)}\}.$
\end{remark}\newpage
\section{A simple recursion for sequence (4)}
Here we prove the following recursive formulas.
\begin{theorem}\label{4}
\begin{equation}\label{15}
x_{2n-1}^{(3)}=2x_{n-1}^{(3)}+1,
\end{equation}
\begin{equation}\label{16}
x_{2n}^{(3)}=2x_{n}^{(3)}.
\end{equation}
\end{theorem}
The initial condition for these recursions are:
\begin{equation}\label{17}
x_1^{(3)}=3, \enskip x_2^{(3)}=6.
\end{equation}
\bfseries Proof. \enskip\mdseries
Let, for a sequence $\{y_n\},$ the formulas (15)-(17) are satisfied:
 \begin{equation}\label{18}
y_{2n-1}=2y_{n-1}+1,
\end{equation}
\begin{equation}\label{19}
y_{2n}=2y_{n},
\end{equation}
\begin{equation}\label{20}
y_1=3, \enskip y_2=6.
\end{equation}
We prove by induction that
\begin{equation}\label{21}
y_n=x_{n}^{(3)}.
\end{equation}
Suppose that (21) is valid for all $n\leq2^t.$ In particular, we have
\begin{equation}\label{22}
y_{2^{t-1}+j}=x_{2^{t-1}+j}^{(3)},\enskip j=0,1,...,2^{t-1}-1
\end{equation}
and (using Lemma 1)
\begin{equation}\label{23}
y_{2^{t}}=x_{2^{t}}^{(3)}=3\cdot2^t.
\end{equation}
It is easy to see that in order to prove (21) it is sufficient to prove (22)-(23) for every $t$ and $j,\enskip j=0,1,...,2^t-1.$ Distinguish cases of odd and even $j.$\newline
$1)$ Let $j$ be odd: $j=1,3,...,2^{t}-1. $ Put in (18) $n=2^{t-1}+\frac {j+1} {2}.$ Using (22) and Lemma 1, we have
$$y_{2^t+j}=2y_{2^{t-1}+\frac {j-1} {2}}+1=$$
$$ 2x_{2^{t-1}+\frac {j-1} {2}}^{(3)}+1=2(3\cdot2^{t-1}+\frac {j-1} {2})+1=3\cdot2^t+j=x_{(2^t+j)}^{(3)},\enskip j=1,3,...,2^{t}-1.$$
$2)$ Let $j$ be even: $j=0,2,...,2^{t}-2. $ Put in (19) $n=2^{t-1}+\frac {j} {2}.$ Using now(23) and Lemma 1, we find
$$y_{2^t+j}=2y_{2^{t-1}+\frac {j} {2}}=$$
$$ 2x_{2^{t-1}+\frac {j} {2}}^{(3)}=2(3\cdot2^{t-1}+\frac {j} {2})=3\cdot2^t+j=x_{(2^t+j)}^{(3)},\enskip j=0,2,...,2^{t}-2.$$
Finally, putting in (19) $n=2^t,$ we have
$$y_{2^{t+1}}=2y_{2^{t}}=2x_{2^{t}}^{(3)}=2(3\cdot2^{t})=3\cdot2^{t+1}=x_{(2^{t+})}^{(3)}. \blacksquare $$
Now we prove an interesting property of sequence (4).\newpage
\begin{corollary}\label{2}
The binary expansion of every number $x_{n}^{(3)}$ does not begin from $10.$
\end{corollary}
\bfseries Proof. \enskip\mdseries The sequence with terms having no $10$ in the beginning of its binary expansions is the Sloane's sequence $A004760 $ [2]. Denote $a(n)=A004760(n), n\geq0.$ In 2004, P. Deleham proved the recursion(see comment to $A004760$):
$$a(0) = 0,\enskip a(1) = 1;$$
for\enskip $n\geq1,$
\begin{equation}\label{24}
 a(2n) = 2a(n) + 1, \enskip a(2n+1) = 2a(n+1) .
\end{equation}
The comparison of (24) with (15)-(17) shows that, for $n\geq1,$ we have: $x_n^{(3)}=a(n+1),$ and the corollary follows. $\blacksquare$\newline
It is interesting that sequence (4) is connected with the very known Josephus combinatorial problem (cf [4]). The solutions of this problem are given by the sequence $A006257.$[2] Denote $b(n)=A006257(n),\enskip n\geq0.$ In 2003, R. Stephan proved that (see comment to $A004760$ [2]): for $n\geq1, \enskip a(n) = 3n - 2 - b(n-1).$ Since $x_n^{(3)}=a(n+1),$  then we have:
\begin{corollary}\label{3}
For $ n\geq1,$ we have
\begin{equation}\label{25}
x_{n}^{(3)}=3n+1-b(n).
\end{equation}
\end{corollary}

\section{Proof of Theorem 3}
First of all, we give an explicit expression for sequence $\{x_n^{(4)}\}$ (see (6)). Here we again denote the exponent of the highest power of 2 dividing $n$ by $(n)_2.$
\begin{lemma}\label{4}
 \begin{equation}\label{26}
x_n^{(4)}=\begin{cases} 2n+3,\;\;if\;\;n((n+1)_2)\enskip is \enskip even,\\2n+2, \;\; otherwise\end{cases}.
\end{equation}
 \end{lemma}
\bfseries Proof. \enskip\mdseries It is easy to see that the statement for $n>1$ is equivalent to the following:
 1) if $n$ is odd, such that either
  $$1a)\enskip (n-1)_2\geq2$$
   or $$1b)\enskip (n-1)_2=1 \enskip and\enskip the \enskip last\enskip series \enskip of \enskip ones\enskip in \enskip the\enskip binary\enskip expression \enskip of\enskip$$ $$n-1\enskip (the\enskip last\enskip1-series)\enskip contains \enskip even\enskip number\enskip of\enskip elements,$$ then $x_n^{(4)}=2n+2;$\newline
 2) if $n$ is odd, such that $ (n-1)_2=1 $ and the last 1-series of $ n-1 $ contains odd number of elements, then $x_n^{(4)}=2n+3;$\newpage
 3)if $n$ is even, then $x_n^{(4)}=2n+3.$\newline
  We prove this modified statement by induction.  For $n=2,3,4,5,6,7,$ where all cases are presented, the formula is true. Denote $ t_n$ the $nth$ Prouhet-Thue-Morse  number [2], i.e. $t_n=0,$ if the number of ones in the binary expansion of $n$ is even, and $t_n=1,$ otherwise. Let the statement be valid for some $n.$\newline
  $a)$Let $n$ be even, such that $(n)_2\geq2.$ Then, by the supposition, $x_n^{(4)}=2n+3.$ Since $(2n)_2\geq3,$ then
  $t(2n+4)=t(2n)+1=t(n)+1=t(n+1) \pmod 2$ and $x_{n+1}^{(4)}=2n+4=2(n+1)+2$ and the lemma follows in subcase $1a).$\newline
$b)$Let $n$ be even, such that $(n)_2=1$ and the last 1-series of $n$ contains even number of ones. Here  $(n)_2=2.$
and we see that, as in $a),$ we have $t(2n+4)=t(2n)+1=t(n)+1=t(n+1) \pmod 2.$ Thus $x_{n+1}^{(4)}=2n+4=2(n+1)+2$ and the lemma follows in subcase $1b).$\newline
$c)$ Let $n$ be even, such that $(n)_2=1$ and the last 1-series of $n$ contains odd number of ones. Here  $(n)_2=2,$
and in this case we ,evidently, have: $t(2n+4)=t(2n)=t(n)\neq t(n+1) \pmod 2,$ but $t(2n+5)=t(2n)+1=t(n)+1=t(n+1) \pmod 2.$ Therefore, $x_{n+1}^{(4)}=2n+5=2(n+1)+3$ and the lemma follows in case 2.\newline
To prove case 3, we distinguish the following subcases:
$d)$ Let $n$ be odd, such that the last 1-series of $n$ contains even number of ones. Here $(n-1)_2=1$ and the last 1-series of $n-1$ contains odd number of ones. Therefore, by the supposition, $x_n^{(4)}=2n+3.$ We have: $t(2n+4)=t(n+2)=t(n)\neq t(n+1).$ On the other hand, $t(2n+5)=t(n)+1=t(n+1) \pmod 2.$ Therefore, $x_{n+1}^{(4)}=2n+5=2(n+1)+3$ and the lemma follows in this subcase of case 3.\newline
$e)$ Let $n$ be odd, such that the last 1-series of $n$ contains odd number of ones. Then, evidently, $t(n+1)=t(n).$ \newline
$e_1)$ The last 1-series of $n$ contains more than 1 ones. Here $(n-1)_2=1$ and the last 1-series of $n-1$ contains even number of ones. Therefore, by the supposition, $x_n^{(4)}=2n+2.$ We have: $t(2n+3)=t(n)+1\neq t(n+1) \pmod 2;$
analogously, $t(2n+4)=t(n+2)=t(n)+1\neq t(n+1) \pmod 2.$ On the other hand, $t(2n+5)=t(n)= t(n+1) $ Therefore, $x_{n+1}^{(4)}=2(n+1)+3$ and the lemma follows in this subcase of case 3.\newline
$e_2)$ The last 1-series of $n$ consists of one 1. Then $(n-1)_2\geq2.$ Therefore, by the supposition, $x_n^{(4)}=2n+2.$ Here $t(2n+3)=t(2n+4)=t(n)+1\neq t(n+1) \pmod 2,$ while $t(2n+5)=t(n)= t(n+1).$ Therefore, $x_{n+1}^{(4)}=2(n+1)+3.$ This completes the proof. $\blacksquare$
\begin{corollary}
 If $n\equiv1 \pmod4,$ then $x_{n}^{(4)}\equiv4 \pmod8,$ and
 \begin{equation}\label{27}
 x_{n+8}^{(4)}-x_{n}^{(4)}=16.
 \end{equation}
\end{corollary}\newpage
\bfseries Proof. \enskip\mdseries By Lemma 4, $x_{n}^{(4)}=2n+2,$ and the statements follow immediately.$\blacksquare$
We shall complete the proof of Theorem 3 in Section 7.
\section{Research of minimal recursive function on 9 consecutive integers}
For some integer $k\geq1,$ consider 9 consecutive integers of the segment $[4k+1,4k+9]$. Let $N\geq1$ be an integer. We
introduce the following integer-valued function $\psi(n)=\psi_{k,N}(n):$ put $\psi(4k+1)=N$ and , if $4k+1<n\leq 4k+9,$ then we consecutively obtain its values by the minimal recursive algorithm respectively property $A_2.$ We want to prove that always $\psi(4k+9)\leq N+16.$ The difficulty consists of the existence of $k, N,$ such that $\psi(4k+5)= N+9.$
\begin{example}\label{5}
Let $k=23,\enskip N=112.$
\end{example}
Then $$\psi(93)=112,\enskip\psi(94)=115,\enskip\psi(95)=116,\enskip\psi(96)=119,\enskip\psi(97)=121.$$
Now we research the possible orders of changes and not-changes of parity of the numbers of binary ones (OCP) of 9 consecutive integers belonging to a $[4k+1,4k+9].$ Denoting every not-change by 0 and every change by 1 \newline( this corresponds to the values of $t_n+t_{n+1} \pmod2$), e.g., in example 5 we have the OCP: $\{0,1,0,1\}.$
\begin{lemma}\label{5}
There are only two OCP of 9 consecutive integers of the segment $[4k+1,4k+9]:\enskip \{0,1,1,1,0,1,1,1\}$ and $\{0,1,1,1.0,1,0,1\}.$
 \end{lemma}
 \bfseries Proof. \enskip\mdseries The last series of ones of $4k+1$ contains, evidently, one 1. It is easy to see that sufficiently to consider series of 0's between the penultimate series of 1's and the last 1, containing 1,2 or 3 zeros, and to fix a parity of 1's in their penultimate series. Further the proof is realized directly by the adding consecutively 1.$\blacksquare$\newline
 It is known (see, e.g., comment by J. O. Shallit to sequence A000069 [2] ("Odious numbers")) that exactly 2 of the 4 numbers  $4t, 4t+1, 4t+2, 4t+3$ have an even sum of binary 1's, while the other 2 have an odd sum. Therefore, the
 change (not-change) of the parity of the number of binary ones always attains by adding of 1,2 or 3 to any integer $n.$
 \begin{definition}
 We call integer $n$ a regular with respect to change (not-change) of the parity of the number of binary ones, if the change (not-change) attains by adding of 1 or 2 to $n.$ Otherwise, $n$ is called a singular with respect to change or not-change correspondingly.\newpage
 \end{definition}
 \begin{lemma}\label{6}
 Every positive integer is regular with respect to change of the parity of the number of binary ones.
 \end{lemma}
\bfseries Proof. \enskip\mdseries If an integer is even, then the statement is trivial. Let an integer be odd with the last series of $m$ 1's. If $m$ is even, then the change of the parity attains by the adding of 1; if $m$ is odd, then the change of the parity attains by the adding of 2.$\blacksquare$\newline
 \begin{lemma}\label{7}
 Every odd positive integer is regular with respect to not-change of the parity of the number of binary ones.
 \end{lemma}
\bfseries Proof. \enskip\mdseries  Let an odd integer have the last series of $m$ 1's. If $m$ is even, then the not-change of the parity attains by the adding of 2; if $m$ is odd, then the not-change of the parity attains by the adding of 1.$\blacksquare$\newline
\begin{lemma}\label{8}
 1)Every even positive integer multiple of 4 is singular with respect  to not-change of the parity of the number of binary ones;\newline
 2) An even positive integer not multiple of 4 is singular with respect to not-change if and only if its last series of 1's has even number of ones.
 \end{lemma}
\bfseries Proof. \enskip\mdseries  Quite analogously.$\blacksquare$
\begin{theorem}\label{5}
For every $k,N\in \mathbb N$ we have $\psi(4k+9)\leq N+16.$
\end{theorem}
\bfseries Proof. \enskip\mdseries We show how the possible "large" jumps of function $\psi(n)$ of the magnitude 3 are compensated by "small" jumps of the magnitude 1. Note that the jumps of function $\psi(n)$ of the magnitude 3 could appear only in 3 points which correspond to 0's of possible OCP according to Lemma 5. Indeed, in other points, by Lemma 6 all integers are regular, therefore, only jumps of the magnitude 1,2 are possible.
1) First OCP. Here the jumps of function $\psi(n)$ of the magnitude 3 appear only in 2 points: $4k+1$ and $4k+5.$ Consider a possibility of appearing of "non-compensating" configuration of jumps of the form $\{3,2,2,2\}.$ on the first singular point. In case of the first type of singulary $...00$ of number $N,$ the first 3 jumps $\{3,2,2,\}$ are possible only in case when the last series of 1's of  $N$ contains odd ones. As a result we obtain a number of the form $...10...011$ with odd last series of 0's. Here, by the first OCP, we have the following change of the parity which, evidently, attains by the adding of 1. Thus,"non-compensating" configuration of jumps  $\{3,2,2,2\}$ is impossible.
In case of the second type of singulary $...01...10$ of number $N$ (here,by Lemma 8, the series of 1 is even)\newpage already after the first 2 jumps $\{3,2\}$ we again obtain a number of the form  $...10...011$ with odd last series of 0's.
Note, that the case of point $4k+5$ is the same. Thus we conclude that the theorem is true in case of the first OCP.
2) Second OCP. Values of function $\psi$ on the segment $[4k+1,4k+5]$ are analyzed by the same way. Thus, we consider
the only OCP for the segment $[4k+5,4k+9]:\enskip \{0,1,0,1\}.$ Here we should consider 4 potential "non-compensating" configurations of jumps  $a)\enskip\{3,2,3,x\},$ where $x=1 or 2,$ $b)\enskip\{3,1,3,2\}$ and $c)\enskip\{3,2,2,2\}.$\newline
$a)$ Independently on a type of singularity, after two first jumps $\{3,2\}$ we obtain an odd number which is always
regular by Lemmas 6-7, therefore, a configuration of jumps  $a)$ is impossible;
$b)$ In case of the first type of singularity of the form $...100$ two first jumps $\{3,1\}$ appear only in case of even last series of 1's and after the first 3 jumps $\{3,1,3\}$ we obtain a number of the form $...10...011$ with positive even number of 0's in the last series of zeros. Here, according to OCP, the following jump is 1. In case of the first type of singularity of the form $...1000$ after the first 3 jumps $\{3,1,3\}$ we obtain a number of the form $...01111$. Here, according to OCP, again the following jump is 1. Finally, in case of the first type of singularity of the form $01...10000$ we indeed obtain a "non-compensating" configuration of jumps $\{3,1,3,2\},$ after which we obtain a number of the form $01...101111.$ Here we use retro-analysis. Subtracting the maximal sum 8 ( as was proved in the above), we obtain an odd number which, by Lemmas 6,7 cannot be singular.Thus, in this case the total sum of jumps is
not more than 16. Consider now the second type of singularity of the form $...01...10,$ where, by Lemma 8, the last series of 1's contains positive even number of ones. After the first jump 3 we obtain a number of the form $10...01$
and the follow jump, according to OCP, cannot be 1. Thus the case $b)$ here is impossible.
$c)$ In case of the first type of singularity of the form $...100$ two first jumps $\{3,2\}$ appear only in case of odd last series of 1's and after them we obtain a number of the form $...10...001$ with positive even number of 0's in the last series of zeros. Now, according to OCP we should add 1. Thus case $c)$ here is impossible. Furthermore, considering the first type of singularity of the form $...1000,$ we see that after the first jump 3 it should be add 1
and case $c)$ here is impossible as well. Finally consider the second type of singularity of the form $...01...10,$ where the last series of 1's contains positive even number of ones. Here we indeed obtain a "non-compensating" configuration of jumps $\{3,2,2,2\},$ after which we obtain a number of the form $...10...0111.$ Using the retro-analysis, we subtract from it the maximal possible sum 8. But we obtain an odd number which cannot be singular.
Thus in this case the total sum of jumps is not more than 16. $\blacksquare$
\begin{corollary}\label{5}
For every $n\geq1,$ we have $x_n^{(2)}<x_n^{(4)}.$
\end{corollary}
\bfseries Proof. \enskip\mdseries  Considering these sequences on positive integers of the form $n=8k+1,$ according
to Theorem 5 and Corollary 2, we conclude that the inequality is true for such $n.$ Now it is sufficient to notice that $x_{33}^{(2)}=51,$ while $x_{33}^{(4)}=68. $ $\blacksquare$
\section{Completion of proof of Theorem 3 }
It is well known that the Prouhet-Thue-Morse sequence is not periodic (a very attractive proof of this fact is given in [3]). We prove a very close statement.
\begin{lemma}\label{9}
There is no a constant $C,$ such that, for every positive integer $n,$ we have $t(2n+C)=t(2n).$
\end{lemma}
 \bfseries Proof. \enskip\mdseries
 Let us take the contrary. Then if $C$ is even, then $C/2$ is a period, which is impossible. If $C$ is odd, i e. $C=2C_1+1,$ then we have
 \begin{equation}\label{28}
1-t(n+C_1)=t(n).
\end{equation}
If $n=2^m,$ where $m>C_1,$ then $t(C_1+n)=1-t(C_1).$ Thus, by (28),we have
 \begin{equation}\label{29}
t(C_1)=1.
\end{equation}
Therefore, if the binary expansion of $C_1$ has the form
 \begin{equation}\label{30}
C_1=2^{r_k}+2^{r_{k-1}}+...+2^{r_1},\enskip r_k>r_{k-1}>...>r_1,
\end{equation}
then $k$ is odd. Consider now $n=2^{r_k}+2^{r_{k}+2}+...+2^{r_{k}+k}.$ Then, the binary expansion of $n$ is $k,$ and, by (30), we have
$$n+C_1=2^{r_{k}+k}+2^{r_{k}+k-1}+...+2^{r_{k}+2}+2^{r_{k}+1}+2^{r_{k-1}}+...+2^{r_1},$$
i.e. the number of 1's in the binary expansion of $n+C_1$ is $2k-1$. Thus, $t(n)=t(n+C_1).$ From this and (28) we have
2t(n)=1. Contradiction. $\blacksquare$\newline
 Moreover, it is easy to see that the equality $t(2n+C)=t(2n)$ cannot be true for every $n\geq n_0.$\newline
Let now $a>4,$ such that $t(a)=1.$  Consider positive integers of the form $n=8k+1.$ By Corollary 4 and Theorem 5, the difference $r(n)=x_n^{(a)}-x_n^{(4)}$ cannot increase. Let us show that it also cannot be constant for $n\geq n_0.$
Indeed, if $r(n)=C,$ then $t( x_n^{(4)}+C)=t(x_n^{(4)}).$  Note that, for the considered form of $n,$ according to Lemma 4, we have $x_n^{(4)}=2n+2.$ Therefore, it should be $t(2(n+1)+C)=t(2(n+1)),$ and, by Lemma 9, it is impossible. Thus, at some moment $r(n)$ attains of the magnitude 1 or 2. Indeed, since the maximal jump of $x_n^{(4)}$ is 3 while the minimal one is 1, then $r(n)$  could \newpage change by jumps 1 or 2. It is left to show that if $r(n)=1,$  then the jump of $r(n)$  could not be 2. Indeed,  since $t(x_n^{(4)})=t(x_n^{(a)}),$ then the case when $x_n^{(4)}$ has jump 3, while   $x_n^{(a)}=x_n^{(4)}+1 $ has jump 1, is impossible, since in the contrary the jump 3 for $x_n^{(4)}$ is not minimal. $\blacksquare$

\;\;\;\;\;\;\;


\begin{thebibliography}{4}
\bibitem 1.  J.-P. \enskip Allouche and J.\enskip Shallit, \slshape The ubiquitous Prouhet-Thue-Morse sequence, \upshape http://www.lri.fr/~allouche/bibliorecente.html
\bibitem 2.  N.\enskip J.\enskip A.\enskip Sloane,\enskip\slshape The On-Line Encyclopedia of Integer Sequences \upshape(http: //www.research.att.com)
\bibitem 3.  S.\enskip  Tabachnikov,  \slshape Variations on Escher theme, \upshape Kvant, no.12 (1990), 2-7 (in Russian).
\bibitem 4.  E.\enskip W.\enskip Weisstein,\enskip\slshape  "Josephus Problem"  From MathWorld--A Wolfram Web Resource (http://mathworld.wolfram.com/JosephusProblem.html)



\end{thebibliography}
\end{document}